\documentclass[12pt,twoside]{article}

\usepackage[utf8]{inputenc}
\usepackage[T1]{fontenc}
\usepackage{amsmath,amssymb,amsfonts,amsthm,mathtools}
\usepackage{enumitem}
\usepackage{microtype}

\allowdisplaybreaks

\newcommand{\C}{\mathbb C}
\newcommand{\N}{\mathbb N}
\newcommand{\R}{\mathbb R}

\theoremstyle{plain}
\newtheorem{theorem}{Theorem}[section]
\newtheorem{lemma}[theorem]{Lemma}
\newtheorem{example}[theorem]{Example}

\newtheorem{corollary}[theorem]{Corollary}
\newtheorem*{theoremA}{Theorem A}
\newtheorem*{theoremB}{Theorem B}
\theoremstyle{definition}

\newtheorem{definition}[theorem]{Definition}
\theoremstyle{remark}

\DeclareMathOperator{\logdens}{logdens}

\numberwithin{equation}{section}

\begin{document}

\title{Complex Oscillation and Non-Oscillation Results of
Second-Order Linear Differential Equations}
\author{Zhi-Tao Wen\footnote{Wen is supported by the National Natural Science Foundation of China (No.~12471076)} and Shi-Jian Wu \footnote{Wu is the corresponding author}}
\date{}
\maketitle

\begin{abstract}
In this paper, we investigate oscillation and non-oscillation results of second order differential equations
    \begin{equation}
    f''+(H_1e^z+H_2e^{\alpha z}+q)f=0,\tag{+}
    \end{equation}
where $H_1,H_2\in\C\setminus\{0\}$, $0<\alpha<1$, and $q$ is an entire function of order less than $1$. First, we show that if equation (+) has a zero-free non-trivial entire solution, then either $\alpha=1/2$ or $\alpha=3/4$, which settles the Problem 9 posed in \cite{HITW}. Second, we prove that if either $q$ is transcendental, or $\alpha=(2m-1)/(2m)$ for some integer $m\geq 2$ and $q\equiv 0$, then every non-trivial entire solution $f$ of $(+)$ satisfies $\lambda(f)=\infty$.
\end{abstract}

\medskip
\noindent\textbf{Keywords.} Complex oscillation; exponential
polynomial; non-oscillation; zero-free solution; exponent of convergence.

\smallskip
\noindent\textbf{2020 MSC.} 34M10, 30D35, 34M05.

\section{Introduction}

Let $A$ be an entire function and consider
\begin{equation}\label{eq:general}
 f''+A(z)f=0.
\end{equation}
It is known that every solution $f$ of \eqref{eq:general} is entire. Moreover, every entire solution $f$ is of finite order if and only if $A$ is a polynomial, and every non-trivial solution $f$ is of infinite order when $A$ is transcendental.
For an entire function $g$, let $\rho(g)$ denote its growth order and
$\lambda(g)$ the exponent of convergence of its zeros.

Heittokangas et al. in \cite{HILT2019}
emphasized‌ the classical oscillation theorem on second
order differential equations \eqref{eq:general} are formulated as:
\begin{quote}
\emph{
Under what conditions on transcendental $A$ does
\eqref{eq:general} have a non-trivial solution $f$ such that
\[
 \textup{[Q1]}.~f~\text{has no zeros},\qquad
 \textup{[Q2]}.~\lambda(f)\ge\rho(A),\qquad
 \textup{[Q3]}.~\lambda(f)=\infty?
\]}
\end{quote}
Bank and Laine \cite{BankLaine1982} proved that for any two linearly independent solutions $f_1$ and $f_2$ of \eqref{eq:general},
the inequality
\[
\max\{\lambda(f_1),\lambda(f_2)\}\geq \rho(A)
\]
holds whenever $\rho(A)<\infty$ is non-integral, and
    $
    \max\{\lambda(f_1),\lambda(f_2)\}=\infty
    $
is true when $\rho(A)<1/2$. In addition, they conjectured that
\[
\max\{\lambda(f_1),\lambda(f_2)\}=\infty
\]
whenever $\rho(A)$ is non-integral, which is known as the
\textbf{Bank--Laine conjecture}. It was later disproved by Bergweiler and
Eremenko in \cite{BergweilerEremenko2017}.

It is much more interesting to study the classical oscillation results
 when $A$ is a specific entire function.
A special case, frequently studied, is given by equations of type
\[
 f''+\bigl(Re^{P}+Q\bigr)f=0,
\]
where $R, P$ and $Q$ are entire functions, see results on \textup{[Q1]}-\textup{[Q3]} in \cite{BankLaineLangley1986}, \cite{BankLangley1987},
\cite{BankLaineLangley1989}, \cite{ChiangLaineWang1995}, \cite{Wu}.

Another special case in \eqref{eq:general} was studied when $A=e^{P_1}+e^{P_2}+Q$, i.e.,
\begin{equation}\label{eq:twoexp}
 f''+\bigl(e^{P_1(z)}+e^{P_2(z)}+Q(z)\bigr)f=0,
\end{equation}
where $P_1$ and $P_2$ are nonconstant polynomials such that $e^{P_1}$ and $e^{P_2}$ are linearly independent over $\C$, and
    $$
    P_1(z)=a_nz^n+\cdots+a_1z+a_0\quad\text{and}\quad P_2(z)=b_nz^n+\cdots+b_1z+b_0,
    $$
and $Q$ is an entire function of order $\rho(Q)<n$. Denote
$\alpha:=a_n/b_n$. The results in
\cite[Theorem~3.2]{IshizakiTohge1997},
\cite[Theorem~1]{Ishizaki1997},
\cite[Theorems~2.1]{HILT2019} and \cite[Corollary 3]{HILT2021} may be summarized as follows.

\begin{theoremA}
Let $f\not\equiv0$ be a solution of \eqref{eq:twoexp}, and let $e^{P_1}$ and $e^{P_2}$ be linearly independent over $\C$. The following
assertions hold.
\begin{enumerate}[label=\textup{\Alph*.}]
\item If $\alpha=1$, then $\lambda(f)\ge n$.
\item If $\alpha<0$ or $\alpha$ is not real, then $\lambda(f)=\infty$.
\item If $0<\alpha<1/2$, then $\lambda(f)\ge n$.
\item If $3/4<\alpha<1$ and $Q\equiv0$, then $\lambda(f)\ge n$.
\item If $P_1=\alpha P_2$, $Q$ is a polynomial with
$\deg Q+2<2n$ for $n\ge2$ and $Q\equiv0$ for $n=1$, and
\[
 0<\alpha<\frac12
 \quad\text{or}\quad
 \frac{2m-1}{2m}<\alpha<\frac{2m+1}{2(m+1)}
\]
for some $m\in\N$, then $\lambda(f)=\infty$.
\end{enumerate}
\end{theoremA}

The results in Theorem A show geometrically that if the points
$0,a_n,b_n$ are collinear, and that $a_n,b_n$ are on the same
side of the origin, then \textup{[Q1]} and \textup{[Q2]} may happen.
The examples in \cite{Ishizaki1997} show that equations
\eqref{eq:twoexp} solve zero-free solutions when $\alpha=1/2$ or $\alpha=3/4$.
Moreover, Ishizaki in \cite{Ishizaki1997} posed a question that is it possible we
can replace $\lambda(f)\geq n$ with $\lambda(f)=\infty$ in the conclusions C and D of
Theorem A.

It inspires us to investigate whether \eqref{eq:twoexp} could have a zero-free solution if
$\alpha$ satisfies
    \begin{equation}\label{m.eq}
    \frac{2m-1}{2m}<\alpha<\frac{2m+1}{2(m+1)}
    \end{equation}
for some $m\in\N$ and $Q\not\equiv 0$. Heittokangas et al. in \cite{HILT2019} focus on this problem and consider
the case $n=1$ which is related to \textup{[Q1]}.

\begin{theoremB}
Suppose that $f$ is a non-constant solution of
\begin{equation}\label{eq:main}
 f''+\left(H_1e^z+H_2e^{\alpha z}+q\right)f=0,
\end{equation}
where $H_1,H_2\in\C\setminus\{0\}$ and $q$ is an entire function of order $\rho(q)<1$.
If $m=1$ or $m=2$ in \eqref{m.eq}, then $f$ is not zero-free.
\end{theoremB}

Following the lengthy proof of Theorem B in \cite{HILT2019}, Heittokangas et al.
mentioned that this result might be extended to the general case.
    \begin{quote}
    \emph{It seems that in the general case $m\in\N$ in \eqref{m.eq} there are no zero-free solutions of \eqref{eq:main} either,
but the method used in this paper can give only one specific interval at the time,
and the amount of considerations increases exponentially along with $m$}.
    \end{quote}
It leads us to investigate whether \eqref{eq:main} admits zero-free solutions for any
$m\in\N$ with $\alpha$ satisfying \eqref{m.eq}. This is stated as Problem 9 in the survey paper \cite{HITW}, which is closely related to \textup{[Q1]}.

\bigskip
\noindent
\textbf{Problem}:
For $m\in\N$, prove that \eqref{eq:main} has no zero-free solution if
$\alpha$ satisfies \eqref{m.eq}.
Find the general forms of the zero-free solutions at
$$
\alpha=\frac{2m-1}{2m},\qquad m\ge3.
$$

\bigskip

The main objective of this paper is to investigate the exponent of convergence of zeros of non-trivial entire solutions of \eqref{eq:main}, and to settle Problem 9 proposed in \cite{HITW}. We state our theorem as follows.

\begin{theorem}\label{thm:classification}
If equation \eqref{eq:main} has a zero-free non-trivial entire solution, then
$\alpha=1/2$ or $\alpha=3/4$. Moreover,
\begin{itemize}
\item If $\alpha=\frac{1}{2}$, and $f$ is a zero-free non-trivial entire solution of
\eqref{eq:main}, then
        $$
        H_1=-a^2,\qquad
        H_2=-a\left(2c+\frac12\right),\qquad
        q\equiv-c^2,
        $$
where $a,c\in\mathbb C$ such that $a\ne0$ and $c\ne-1/4$. The zero-free solution is
    $$
    f=K\exp\left(cz+2ae^{\frac{z}{2}}\right),
    $$
where $K\in\C\setminus\{0\}$.
\item If $\alpha=\frac{3}{4}$, and $f$ is a zero-free non-trivial entire solution of
\eqref{eq:main}, then
        $$
        H_1=-16c^4,\qquad
        H_2=8c^3,\qquad
        q\equiv-\frac1{64},
        $$
where $c\in\mathbb C\setminus\{0\}$. The zero-free solution is
    $$
    f=K\exp\left(-\frac{z}{8}-8c^2e^{\frac{z}{2}}+4ce^{\frac{z}{4}}\right),
    $$
where $K\in\C\setminus\{0\}$.
\end{itemize}
\end{theorem}

A  further question related to \textup{[Q3]} is
under what condition on $A$, every non-trivial entire solution of \eqref{eq:main} satisfies $\lambda(f)=\infty$. The following result shows that for almost cases, every non-trivial entire solution of \eqref{eq:main} satisfies $\lambda(f)=\infty$.

\begin{theorem}\label{thm:infinity2}
Suppose that $q$ is a transcendental entire function of order $\rho(q)<1$, and $\alpha\in\R$ such that $0<\alpha<1$. Then every non-trivial solution $f$ of \eqref{eq:main} satisfies $\lambda(f)=\infty$.
\end{theorem}

According to the statement in \cite[Page~1981]{Zhang}, Zhang proved that
all non-trivial solutions of equation \eqref{eq:main} satisfy $\lambda(f)=\infty$, when $q$ is a non-constant polynomial and $\alpha\in\mathbb{Q}$ such that $0<\alpha<1$ by using the method in \cite{Langley}. Together with Theorem \ref{thm:infinity2}, we have the following corollary.

\begin{corollary}
Suppose that $q$ is a non-constant entire function of order $\rho(q)<1$, and $\alpha\in\mathbb{Q}$ such that $0<\alpha<1$. Then every non-trivial solution $f$ of \eqref{eq:main} satisfies $\lambda(f)=\infty$.
\end{corollary}

Note that \textup{[Q1]}, \textup{[Q2]} and \textup{[Q3]} might happen when $q$ is a constant in \eqref{eq:main}, see, e.g., \cite[Example~1.1]{HILT2019}. In the following, we give a complex oscillation result of \eqref{eq:main} in a special case, which is related to \textup{[Q3]}.

\begin{theorem}\label{thm:infinity}
Let $m\in\mathbb{N}^+$ and $\alpha=\frac{2m-1}{2m}$. Suppose that $m\geq 2$ and $q=0$.
Then every non-trivial solution $f$ of equation \eqref{eq:main} satisfies $\lambda(f)=\infty$.
\end{theorem}

We note here that the assumption $m\geq 2$ in Theorem \ref{thm:infinity} is sharp. For example,
the entire function $f(z)=\exp(e^{2iz})$ satisfies
    $$
    f''+(4e^{4iz}+4e^{2iz})f=0.
    $$
It is obvious that $\lambda(f)=0<\infty$. In addition, it is impossible to remove $q\equiv 0$.
For example, the differential equation
    $$
    f'' - \left(\frac{1}{16} e^{z} + \frac18 e^{3z/4} + \frac{25}{64}\right) f = 0
    $$
has an entire solution
$
f(z) = \bigl(1 - e^{z/4}\bigr) \exp\left(\frac12 e^{z/2} + e^{z/4} - \frac58 z\right)
$
with $\lambda(f)=1<\infty$. We can construct more examples using \cite[Theorem~4.4]{Zhang}.

\section{Proof of Theorem~\ref{thm:classification} }

In order to prove Theorem~\ref{thm:classification}, let us recall a key lemma, which is from \cite[Theorem~2.1]{ZhangCMFT}. We state it as follows.
\begin{lemma}\label{prop:frequency}
    Let $B_1B_2\ne0$, $0<\alpha<1$, and let $q$ be an entire function of order
    $\rho(q)<1$. If an entire function $u$ is a solution of
    $$
     u^2+u'+q=B_1e^z+B_2e^{\alpha z},
    $$
then
    $$
     u(z)=\gamma(z)+\sum_{j=0}^m a_je^{\mu_jz},
    $$
where $\gamma$ is an entire function of order $\rho(\gamma)<1$, $a_0,\ldots,a_m$ are non-zero constants, $\mu_j=\frac{1}{2}-j(1-\alpha)$, $j=0,\ldots,m$, and
    $$
    m=\bigg\lceil\frac{1}{2(1-\alpha)}\bigg\rceil-1,
    $$
where $\lceil a\rceil$ is the smallest integer not less than $a$.
    \end{lemma}
Suppose that $f$ is a zero-free non-trivial entire solution of \eqref{eq:main}. There exists
an entire function $g$ such that $f=e^g$. Denote
    $$
    F:=\frac{f'}f=g'.
    $$
Thus, it yields that
\begin{equation}\label{eq:riccatiF}
 F'+F^2=-H_1e^z-H_2e^{\alpha z}-q.
\end{equation}
Using Lemma~\ref{prop:frequency} and \eqref{eq:riccatiF}, we have
    $$
    F=\gamma(z)+\sum_{j=0}^m a_je^{\mu_jz},
    $$
where $\mu_j=\frac12-j(1-\alpha)$ and $a_j\in\C\setminus\{0\}$, and $\gamma$ is an entire function of order $\rho(\gamma)<1$. Obviously, we see that $0<\mu_m<\cdots<\mu_1<\mu_0=\frac{1}{2}$. In addition,
    \begin{equation}\label{mu.eq}
    \mu_j-\alpha=\frac{1}{2}-j(1-\alpha)-\alpha
                =-\frac{1}{2}+(1-j)(1-\alpha)<0
    \end{equation}
for $1\leq j\leq m$.
\begin{equation}\label{eq:Fexpand}
 F'+F^2=\gamma'+\gamma^2+\sum_{j=0}^m a_j(2\gamma+\mu_j)e^{\mu_jz}+\sum_{j=0}^m\sum_{\ell=0}^m
 a_ja_\ell e^{(\mu_j+\mu_{\ell})z}
\end{equation}
Let us discuss \eqref{eq:Fexpand} in the following three cases.

\bigskip
\textbf{Case 1: $m\geq 2$}. Obviously, we have $\alpha> 3/4$. Since $\mu_m<\mu_j$ for $0\leq j\leq m-1$, and $\mu_m<\mu_j+\mu_{\ell}$ for $0\leq j,\ell\leq m$, using \cite[Theorem~5.1]{YangYi} and \eqref{eq:riccatiF}-\eqref{eq:Fexpand}, we obtain
    $$
    2\gamma+\mu_m=0\quad\text{and}\quad \gamma^2+\gamma'=-q.
    $$
It implies that $\gamma$ and $q$ are constants. Thus, the coefficients of $e^{\mu_jz}$ cannot vanish for $j=0,\ldots,m-1$. Using \cite[Theorem~5.1]{YangYi} and \eqref{eq:Fexpand} again, we have
    $$
    \mu_{m-1}=2\mu_m.
    $$
Obviously, it gives that
    $$
    \delta:=1-\alpha=\frac{1}{2(m+1)}.
    $$
Denote $N=m+1$ and $x=e^{\delta z}$.
Hence, we see that
    $$
    e^{\mu_jz}=e^{\left(\frac{1}{2}-j(1-\alpha)\right)z}=e^{(N-j)\delta z}=x^{N-j}
    $$
for $j=0,\ldots,N-1$. We rewrite $F$ as
    \begin{equation}\label{F.eq}
    F(z)=\gamma+\sum_{j=0}^{m}a_je^{\mu_jz}=\gamma+\sum_{j=0}^ma_jx^{N-j}
    =\sum_{j=0}^Nr_jx^j=:R(x),
    \end{equation}
where $r_0=\gamma$ and $r_j=a_{N-j}\neq 0$ for $j=1,\ldots,N$. Substitute \eqref{F.eq} into
\eqref{eq:riccatiF}, we have
    \begin{equation}\label{eq:polyriccati}
    R(x)^2+\delta xR'(x)
    =-H_1x^{2N}-H_2x^{2N-1}-q.
    \end{equation}
It yields by \eqref{F.eq} and \eqref{eq:polyriccati} that
    $$
    r_{N}^2=-H_1,\quad 2r_Nr_{N-1}=-H_2,\quad r_0^2=-q.
    $$
Moreover, denote $c_n=\sum_{j+k=n}r_jr_k$. Hence, it gives that
    \begin{equation}\label{cn1.eq}
    c_n=0 \quad \text{for} \quad N+1\leq n\leq 2N-2,
    \end{equation}
and
    \begin{equation}\label{cn2.eq}
    c_n+\delta r_n n=0 \quad\text{for}\quad 1\leq n\leq N.
    \end{equation}
It follows from \eqref{cn1.eq} by mathematical induction that
    \begin{equation}\label{eq:highcoeff}
    r_{N-k}=\frac{1}{k}\frac{(-1)^{k+1}}{2^{k-1}}\binom{2k-2}{k-1}\frac{r_{N-1}^k}{r_N^{k-1}}
    \end{equation}
for $k=2,\ldots,N-1$. Choosing $n=1,2,3$ in \eqref{cn2.eq}, we have
    \begin{equation}\label{r2}
    r_1^2=2r_0r_2,\quad
    r_2^2=r_1r_3.
    \end{equation}
\textbf{Subcase $N\geq 5$}. In particular, if $k=N-1,N-2,N-3$ in \eqref{eq:highcoeff}, then
    $$
    \frac{r_2^2}{r_1r_3}=\frac{(N-1)(2N-7)}{(N-2)(2N-5)}.
    $$
From \eqref{r2}, we know that $(N-1)(2N-7)=(N-2)(2N-5)$, which is impossible.

\bigskip
\noindent
\textbf{Subcase $N=4$}. From \eqref{cn1.eq}, we see that $c_n=0$ for $5\leq n\leq 6$. It gives that
    $$
    r_1r_4+r_2r_3=0,\quad r_3^2+2r_2r_4=0.
    $$
It follows that $r_2=0$ by \eqref{r2}, which is impossible.

\bigskip
\noindent
\textbf{Subcase $N=3$}. From \eqref{cn1.eq}, we see that $c_4=0$. It gives that
    $$
    2r_1r_3+r_2^2=0.
    $$
It follows that $r_2=0$ by \eqref{r2}, which is impossible.

\bigskip
Therefore, it is a contradiction with $m\geq 2$.

\bigskip
\textbf{Case 2: $m=1$}. It is obvious that $1/2< \alpha\leq 3/4$.
From \eqref{eq:riccatiF} and \eqref{eq:Fexpand}, we get
    $$
    \gamma'+\gamma^2+a_0\left(2\gamma+\frac{1}{2}\right)e^{\frac{z}{2}}+ a_1(2\gamma+\mu_1)e^{\mu_1z}+a_0^2e^{z}+2a_0a_1e^{\alpha z}+
    a_1^2e^{2\mu_1z}=-H_1e^z-H_2e^{\alpha z}-q.
    $$
Using \cite[Theorem~5.1]{YangYi}, we have $\alpha=3/4$, and
    $$
    2\gamma+\mu_1=0, \quad a_0^2=-H_1, \quad a_0\left(2\gamma+\frac{1}{2}\right)+a_1^2=0,\quad
    2a_0a_1=-H_2,\quad
    \gamma'+\gamma^2=-q.
    $$
Therefore, it gives that
        $$
        H_1=-16c^4,\qquad
        H_2=8c^3,\qquad
        q\equiv-\frac1{64},
        $$
where $c\in\mathbb C\setminus\{0\}$. The zero-free solution is
    $$
    f=K\exp\left(-\frac{z}{8}-8c^2e^{\frac{z}{2}}+4ce^{\frac{z}{4}}\right),
    $$
where $K\in\C\setminus\{0\}$.

\bigskip
\textbf{Case 3: $m=0$}. It is obvious that $0<\alpha\leq 1/2$.
From \eqref{eq:riccatiF} and \eqref{eq:Fexpand}, we get that
    $$
    \gamma'+\gamma^2+a_0\left(2\gamma+\frac{1}{2}\right)e^{\frac{z}{2}}+ a_0^2e^{z}=-H_1e^z-H_2e^{\alpha z}-q.
    $$
Using \cite[Theorem~5.1]{YangYi}, we have $\alpha=1/2$, and
    $$
    a_0^2=-H_1, \quad a_0\left(2\gamma+\frac{1}{2}\right)=-H_2,\quad
    \gamma'+\gamma^2=-q.
    $$
Therefore, it gives that
        $$
        H_1=-a^2,\qquad
        H_2=-a\left(2c+\frac12\right),\qquad
        q\equiv-c^2,
        $$
where $a,c\in\mathbb C$ such that $a\ne0$ and $c\ne-1/4$. The zero-free solution is
    $$
    f=K\exp\left(cz+2ae^{\frac{z}{2}}\right),
    $$
where $K\in\C\setminus\{0\}$. This completes the proof.

\section{Proof of Theorem~\ref{thm:infinity2}}
Assume to the contrary that $\lambda(f)<\infty$. Then
$$
f=\Pi e^g,
$$
where $\Pi$ is entire of finite order, and $g$ is an entire function of order 1 according to \cite[Lemma 2.10]{Chen2003} and \cite[Theorem~1]{Won}. Put $u=g'$. Consequently, $\rho(u)=1$.
Substituting into
\eqref{eq:main} gives
\begin{equation}\label{eq:transdirect}
u^2+u'+2\frac{\Pi'}{\Pi}u+\frac{\Pi''}{\Pi}
+H_1e^z+H_2e^{\alpha z}+q=0.
\end{equation}
Using \cite[Corollary~5.13]{Laine},
there exist $K>0$ and an $R$-set $E$ such that for all sufficiently large $z\not\in E$, we have
\begin{equation}\label{eq:translogder}
\left|\frac{\Pi'}{\Pi}\right|
+\left|\frac{\Pi''}{\Pi}\right|
+\left|\frac{u'}u\right|\le |z|^K.
\end{equation}
We note that the following $K>0$ stands for a constant, not necessarily the same at each occurrence.

First, let us consider the asymptotic growth of $u$ in the interval $-\frac{\pi}{2}<\theta<\frac{\pi}{2}$. From \eqref{eq:transdirect} and \eqref{eq:translogder}, we have there exists $K>0$ such that for all sufficiently large $z\not\in E$,
    \begin{equation}\label{H.eq}
    |H_1e^z+H_2e^{\alpha z}+q|\leq
    |u|^2+|u'|+2\left|\frac{\Pi'}{\Pi}u\right|+\left|\frac{\Pi''}{\Pi}\right|
    \leq |u|^2+|z|^K|u|+|z|^K.
    \end{equation}
Assume that there exists a sequence $z_n$ with $|z_n|\to\infty$ such that $|u(z_n)|\leq |z_n|^K$. Then for all sufficiently large $z_n\not\in E$, \eqref{H.eq} gives
    $$
    |H_1e^{z_n}+H_2e^{\alpha z_n}+q(z_n)|=O(|z_n|^K),
    $$
 which is impossible. Therefore, for any constant $N$ and all sufficiently large $z\not\in E$, we have
    $$
    |u(z)|> |z|^N.
    $$
By \eqref{H.eq} again, it yields that all sufficiently large $z\not\in E$
    $$
    |H_1e^{z}+H_2e^{\alpha z}+q|\leq |u|^2(1+o(1)).
    $$
Similarly, for all sufficiently large $z\not\in E$, it gives that
    $$
    |H_1e^{z}+H_2e^{\alpha z}+q|\geq |u|^2(1-o(1)).
    $$
Therefore, we obtain that in the interval $-\frac{\pi}{2}<\arg z=\theta<\frac{\pi}{2}$,
    \begin{equation}\label{u1.eq}
    |u(z)|\asymp e^{\frac{1}{2}|z|\cos\theta}\quad \text{as}~r\to\infty,
    \end{equation}
for all sufficiently large $z\not\in E$.

Next, we claim that in the interval $\frac{\pi}{2}<\arg z=\theta<\frac{3\pi}{2}$,
for any $\varepsilon>0$,
    \begin{equation}\label{u2.eq}
    |u(z)|=O\left(e^{r^{\rho(q)+\varepsilon}}\right),
    \end{equation}
for all sufficiently large $z\not\in E$. Otherwise, there exists a sequence $\{z_n\}$ outside of $E$ with $|z_n|\to\infty$ such that
    $$
    |u(z_n)|> e^{|z_n|^{\rho(q)+\varepsilon}}> \frac{1}{4}|z_n|^K.
    $$
It shows by \eqref{eq:transdirect} and \eqref{eq:translogder} that
    $$
    \frac{1}{4}|u(z_n)|^2< |u(z_n)|^2-|u'(z_n)|-2|z_n|^K|u|-|z_n|^K
    \leq |H_1e^{z_n}+H_2e^{\alpha z_n}+q(z_n)|\leq e^{|z_n|^{\rho(q)+\varepsilon}},
    $$
which is a contradiction. Hence, our claim \eqref{u2.eq} is true.

In the following , we put
\[
\delta=1-\alpha,
\qquad
m=\min\left\{j\ge0:\delta\ge\frac1{2(j+1)}\right\}.
\]
Denote $a^2=-H_1$ and $\mu_j=1/2-j\delta$, and set
\begin{equation}\label{eq:transR}
R(z):=a\sum_{j=0}^{m}\binom{1/2}{j}
\left(\frac{H_2}{H_1}\right)^j
e^{\mu_jz}=\sum_{j=0}^mc_je^{\mu_jz}.
\end{equation}
In addition, let us set
\begin{equation}\label{eq:transE}
\begin{split}
E(z):&=R(z)^2+H_1e^z+H_2e^{\alpha z}
=\sum_{j=0}^m\sum_{k=0}^m c_j c_k e^{(\mu_j+\mu_k)z}+H_1e^z+H_2e^{\alpha z}\\
&=\sum_{n=0}^{2m}\left(\sum_{\substack{j+k=n\\0\le j,k\le m}}c_j c_k\right)e^{(1-n\delta)z}
+H_1e^z+H_2e^{\alpha z}\\
&=\sum_{n=0}^{2m}A_ne^{(1-n\delta)z}
+H_1e^z+H_2e^{\alpha z}.
\end{split}
\end{equation}
Since $A_0=a^2=-H_1$ and $A_1=-H_2$, and
$$
\displaystyle\sum_{j+k=n}\binom{\frac12}{j}\binom{\frac12}{k}=0
$$
for $2\leq n\leq m$, see \cite[Page~84]{Gould}, we have
    $$
    E(z)=\sum_{n=m+1}^{2m}A_ne^{(1-n\delta)z}.
    $$
We note that for $n\geq m+1$
    $$
    1-n\delta\leq 1-(m+1)\delta\leq \frac{1}{2}.
    $$
Using \eqref{eq:transR} and \eqref{eq:transE}, we rewrite \eqref{eq:transdirect} as
\begin{equation}\label{eq:transu}
(u-R)(u+R)+u'+2\frac{\Pi'}{\Pi}u
+\frac{\Pi''}{\Pi}+q+E=0.
\end{equation}
Since $E$ is an $R$-set, the set of angles for which the rays meets
infinitely many discs of an $R$-set has linear measure zero, see \cite[Lemma~5.9]{Laine}. According to $\rho(q)<1$, we may choose $\varepsilon$ so small satisfying
$\pi/2\pm\varepsilon, 3\pi/2\pm\varepsilon$ do not meet $E$ for large $z$, such that
\begin{equation}\label{q.eq}
 \rho(q)+\varepsilon<\frac{\pi}{\pi+2\varepsilon}<1.
\end{equation}
For $|\arg z-\pi|\le\pi/2-\varepsilon$, where $\varepsilon$ is defined in \eqref{q.eq},
we obtain that for all sufficiently large $z$,
\begin{equation}\label{eq:transleft}
\log^+|R(z)|+\log^+|E(z)|+\log^+|q(z)|=O(|z|^{\rho(q)+\varepsilon}).
\end{equation}
For $|\arg z|\le\pi/2-\varepsilon$, where $\varepsilon$ is defined in \eqref{q.eq}, we see that in this sector,
\begin{equation}\label{eq:transright}
|R(z)|=|a|e^{\Re z/2}(1+o(1)),
\qquad |E(z)|+|q(z)|=O(|R(z)|).
\end{equation}
Denote $G:=\{z\not\in E: |\arg z|\le\pi/2-\varepsilon\}$.
Equations \eqref{eq:transu}, \eqref{eq:translogder}, and
\eqref{eq:transright} give, there exists $K>0$ such that
\[
|(u(z)-R(z))(u(z)+R(z))|
\le |z|^K\bigl(|u(z)|+|R(z)|+1\bigr)
\]
for all sufficiently large $z$ satisfying $z\in G$. Since
\[
|u|+|R|\le |u-R|+|u+R|,
\qquad \max\{|u-R|,|u+R|\}\ge |R|,
\]
and division by the larger
factor gives
\[
\min\{|u-R|,|u+R|\}\le3|z|^K,
\]
for all sufficiently large $z$ satisfying $z\in G$.
Thus, there exists $K>0$ such that
\begin{equation}\label{eq:transalternative}
|u(z)-R(z)|\le |z|^K
\quad\text{or}\quad
|u(z)+R(z)|\le |z|^K
\end{equation}
for all sufficiently large $z$ satisfying $z\in G$.
Let us set
    $$
    h(z)=|u-R|-|u+R|.
    $$
Obviously, $h(z)$ is a real continuous function. Moreover, 
from Lemma \ref{connect.lemma}, we know that $G$ is a connected set by choosing a suitable $R$-set $E$.
According to \eqref{eq:transright}, for all $z\in G$, the two inequalities in \eqref{eq:transalternative} cannot hold simultaneously, otherwise we got $2|R|\le |u-R|+|u+R|\le2|z|^K$, which is a contradiction with \eqref{eq:transright}. It shows that $h(z)$ never vanish for all $z\in G$. Denote
    $$
    E_1:=\{z\in G :|u(z)-R(z)|\leq |z|^K\},\quad
    E_2:=\{z\in G :|u(z)+R(z)|\leq |z|^K\}.
    $$
According to the discussion above, we have $G=E_1$ or $G=E_2$.

Suppose that $G=E_1$. From \eqref{eq:translogder}, \eqref{eq:transu} and
\eqref{eq:transleft}, we know that
    $$
    |u-R|=O\left(e^{|z|^{\rho(q)+\varepsilon}}\right)
    $$
for all sufficiently large $z\not\in E$ in the sector $\frac{\pi}{2}\leq \arg z\leq \frac{3z}{2}$.
Since there exists $K>0$ such that
    $$
    |u-R|=O(|z|^K)
    $$
for all sufficiently large $z\not\in E$ in the sector $-\frac{\pi}{2}\leq \arg z\leq \frac{z}{2}$,
By means of the Phragm\'{e}n-Lindel\"{o}f theorem, for any $0\leq \arg z<2\pi$,
we conclude that
    $$
    |u-R|=O\left(e^{|z|^{\rho(q)+\varepsilon}}\right)
    $$
for all sufficiently large $z$.
By means of \eqref{q.eq} and the Phragm\'{e}n-Lindel\"{o}f theorem in the sector, see
\cite[Theorem~4.11]{GH2026}, we found that $u-R$ is a polynomial, because
 there exists $K>0$ such that
    $$
    |u-R|=O(|z|^K)
    $$
for all sufficiently large $z$ along the rays $\arg z=\pi/2-\varepsilon$ and
$\arg z=3\pi/2+\varepsilon$. Hence, $u$ and $u'$ are both polynomial growth in the sector $\frac{\pi}{2}<\arg \pi<\frac{3\pi}{2}$. By \eqref{eq:transdirect} and \eqref{eq:translogder}, it gives us that
there exists $K>0$ such that
    $$
    |q|=O(|z|^K)
    $$
for all sufficiently large $z\not\in E$ satisfying $\frac{\pi}{2}\leq \arg \pi\leq \frac{3\pi}{2}$. According to \eqref{q.eq} and Phragm\'{e}n-Lindel\"{o}f theorem, see
\cite[Theorem~4.11]{GH2026}, we conclude that $q$ is a polynomial. It is a contradiction.

Similarly, we also get a contradiction if $G=E_2$. Therefore, we proved our assertion.

\section{Proof of Theorem~\ref{thm:infinity}}

Suppose that equation \eqref{eq:main} with $q=0$ and
$\alpha=(2m-1)/(2m)$, where $m\in\N$ such that $m\geq 2$, admits a non-trivial solution $f$ such that $\lambda(f)<\infty$. Put
\[
 z=2mw,\qquad y(w)=f(2mw).
\]
Denote $y'=dy/dw$. We have $y$ is a solution of
\begin{equation}\label{y.eq}
y''(w)+4m^2\bigl(H_1e^{2mw}+H_2e^{(2m-1)w}\bigr)y(w)=0.
\end{equation}
It is obvious that the coefficient of the equation \eqref{y.eq} is a non-constant entire
function of period $2\pi i$ and is a polynomial in $e^w$. Hence, we know the representation of $f$ by using \cite[Theorem~1]{BankLaine}.

\begin{lemma}\label{Laine.lemma}
Let $A(z)$ be a nonconstant periodic entire function of period $\omega$ and rational in $e^{\alpha z}$, where $\alpha=2\pi i\omega^{-1}$. Let $f(z)\not\equiv 0$ be a solution of $f''+A(z)f=0$, such that the exponent of convergence of the zero-sequence of $f(z)$ is finite. Then the following hold:

\begin{itemize}
\item
If $f(z)$ and $f(z+\omega)$ are linearly dependent, then
\[
f(z)=\Psi(e^{\alpha z})\exp\left(\sum_{j=q}^{m} d_j e^{\alpha j z}+vz\right),
\]
where $\Psi(\rho)$ is a polynomial with all simple nonzero roots, $m,q\in\mathbb Z$, $m\ge q$, $v,d_q,\dots,d_m\in\mathbb C$, $d_j\neq 0$ for some $j\neq 0$.
\item
If $f(z)$ and $f(z+\omega)$ are linearly independent, then
\[
f(z)=\Psi(e^{\frac{\alpha}{2} z})\exp\left(e^{\frac{\alpha}{2} z}\sum_{j=q}^{m} d_j e^{\alpha j z}+vz\right),
\]
where $\Psi(\rho)$ is a polynomial with all simple nonzero roots, $m,q\in\mathbb Z$, $m\ge q$, $v,d_q,\dots,d_m\in\mathbb C$, $d_j\neq 0$ for some $j$.

\end{itemize}
\end{lemma}

Next, we proceed to prove our assertion in the following two cases.

\textbf{Case 1}. $y(w)$ and $y(w+2\pi i)$ are linearly independent. According to Lemma~\ref{Laine.lemma}, we have that
\[
 y(w)=\Psi\left(e^{\frac{w}{2}}\right)\exp\left(
 e^{\frac{w}{2}}\sum_{j=\ell}^{s}d_je^{jw}+\nu w\right),
\]
where $\Psi(z)$ is a polynomial all of whose roots are simple and nonzero,
$\ell \le s$ are integers, and
$d_{\ell},\ldots,d_s,\nu\in\C$ such that $d_{\ell}d_s\ne0$.
Thus,
\begin{align*}
\frac{y''}{y}
&=\frac{e^{\frac{w}{2}}}{4}\frac{\Psi'(e^{\frac{w}{2}})}{\Psi(e^{\frac{w}{2}})}
+\frac{e^w}{4}\frac{\Psi''(e^{\frac{w}{2}})}{\Psi(e^{\frac{w}{2}})}
+\frac14\sum_{j=\ell}^s d_j(2j+1)^2 (e^{\frac{w}{2}})^{2j+1}\\
&+\frac{e^{\frac{w}{2}}\Psi'(e^{\frac{w}{2}})}{\Psi(e^{\frac{w}{2}})}\left(\frac12\sum_{j=\ell}^s d_j(2j+1)(e^{\frac{w}{2}})^{2j+1}+\nu\right)
+\left(\frac12\sum_{j=\ell}^s d_j(2j+1)(e^{\frac{w}{2}})^{2j+1}+\nu\right)^2.
\end{align*}

If \(s+\tfrac12\le 0\), then \(\sum_{j=\ell}^s d_j(j+\tfrac12)e^{(j+1/2)w}\) remains bounded as \(w\to+\infty\). Consequently, \(\tfrac{y''}{y}\) is bounded as \(w\to+\infty\). This contradicts the asymptotic behaviour \(\tfrac{y''}{y}\sim -4m^2H_1 e^{2mw}\) as \(w\to+\infty\). Therefore we must have \(s+\tfrac12>0\). Since \(d_s\neq 0\), the dominant term as \(w\to+\infty\) is
\(\frac{y''}{y}\sim \left(s+\frac12\right)^2 d_s^2 e^{(2s+1)w}.\)
Comparing with \(\tfrac{y''}{y}\sim -4m^2\left(H_1 e^{2mw}\right)\), we obtain \(2s+1=2m\). This is impossible, since the left-hand side is odd and the right-hand side is even.

\medskip
\textbf{Case 2}. $y(w)$ and $y(w+2\pi i)$ are linearly dependent.
According to Lemma~\ref{Laine.lemma}, we have that
\[
 y(w)=\Psi(e^w)\exp\left(
 \sum_{j=\ell}^{s}d_je^{jw}+\nu w\right),
\]
where $\Psi(z)$ is a polynomial all of whose roots are simple and nonzero,
$\ell \le s$ are integers, and
$d_{\ell},\ldots,d_s,\nu\in\C$ such that $d_{\ell}d_s\ne0$.
Thus,

\begin{align}\nonumber
\frac{y''}{y}
=&e^w\frac{\Psi'(e^w)}{\Psi(e^w)}+e^{2w}\frac{\Psi''(e^w)}{\Psi(e^w)}+\sum_{j=\ell}^s j^2 d_j e^{jw}\\\label{y''.eq}
+&2e^w\frac{\Psi'(e^w)}{\Psi(e^w)}\left(\sum_{j=\ell}^s j d_j e^{jw}+\nu\right)
+\left(\sum_{j=\ell}^s j d_j e^{jw}+\nu\right)^2.
\end{align}

Suppose that \(\ell< 0\). Then \(\tfrac{y''}{y}\) is unbounded as \(w\to-\infty\), which contradicts the asymptotic behaviour \(\tfrac{y''}{y}\sim -4m^2H_2e^{(2m-1) w}\to 0\) as \(w\to-\infty\). Hence we must have \(\ell\geq 0\) and $v^2=0$. Without loss of generality, we assume that $\ell>0$.

Assume that $s=\ell$. It follows that $2s=2m$ and $s=2m-1$ by \eqref{y''.eq}.
The only possibility is $m=1$. It is a contradiction with our assumption.
Hence, it yields that $s-\ell\geq 1$. As \(w\to+\infty\),
$$
\frac{y''}{y}= s^2 d_s^2 e^{2sw}+2s(s-1)d_sd_{s-1}e^{(2s-1)w}+O(e^{(2s-2)w}).
$$
Comparing with equations \eqref{y.eq}, we get \(s=m\) and $d_{s-1}\neq 0$.
Hence,
    $$
    \left(\sum_{j=\ell}^s j d_j e^{jw}\right)^2+4m^2(H_1e^{2mw}+H_2e^{(2m-1)w})=O(e^{mw}),
    $$
as $w\to+\infty$. Using $(1+x)^{1/2}=\sum_{j=0}^{\infty}\binom{\frac12}{j}x^{j}$ for $|x|<1$ and setting $a^2=-H_1$, we obtain that
    \begin{align}\nonumber
    \sum_{j=\ell}^s j d_j e^{jw}&=2mae^{mw}\sum_{n=0}^\infty\dbinom{\frac{1}{2}}{n}\left(\frac{H_2}{H_1}e^{-w}+O(e^{-mw})\right)^n\\\label{sj.eq}
    &=2mae^{mw}\sum_{n=0}^\infty\sum_{k=0}^n\dbinom{\frac{1}{2}}{n}\dbinom{n}{k}
    \left(\frac{H_2}{H_1}e^{-w}\right)^{n-k}\left(O(e^{-mw})\right)^k,
    \end{align}
as $w\to+\infty$. Since $0<\ell\leq m$, comparing the positive powers of $e^w$ in \eqref{sj.eq} yields
    \begin{equation}\label{jdj.eq}
    jd_j=2ma\dbinom{\frac{1}{2}}{m-j}\left(\frac{H_2}{H_1}\right)^{m-j}
    \end{equation}
for $j=\ell,\ldots,m=s$.

Suppose that $s=m=2p$. It is obvious that $p\geq 1$. Using \eqref{jdj.eq}, we have the coefficients of $e^{mw}$ in $\left(\sum_{j=\ell}^m j d_j e^{jw}\right)^2$ is
    $$
    2\sum_{t=1}^{p-1}t(m-t)d_td_{m-t}+(pd_p)^2
    =8m^2a^2\left(\frac{H_2}{H_1}\right)^{m}\left(\sum_{t=1}^{p-1}
    \dbinom{\frac{1}{2}}{m-t}\dbinom{\frac{1}{2}}{t}+\frac{1}{2}\dbinom{\frac{1}{2}}{p}^2\right).
    $$
In addition, \cite[Page~84]{Gould} gives
    $$
    0=\sum_{t=0}^{2p}\dbinom{\frac{1}{2}}{m-t}\dbinom{\frac{1}{2}}{t}
    =2\dbinom{\frac{1}{2}}{2p} +\dbinom{\frac{1}{2}}{p}^2+\sum_{t=1}^{p-1}\dbinom{\frac{1}{2}}{m-t}\dbinom{\frac{1}{2}}{t}
   +\sum_{t=p+1}^{2p-1}\dbinom{\frac{1}{2}}{m-t}\dbinom{\frac{1}{2}}{t}.
    $$
It yields that
    $$
    \sum_{t=1}^{p-1}
    \dbinom{\frac{1}{2}}{m-t}\dbinom{\frac{1}{2}}{t}+\frac{1}{2}\dbinom{\frac{1}{2}}{p}^2
    =-\dbinom{\frac{1}{2}}{2p} =-\dbinom{\frac{1}{2}}{m}.
    $$
Suppose that $s=m=2p+1$. It is obvious that $p\geq 1$. Using \eqref{jdj.eq}, we have the coefficients of $e^{mw}$ in $\left(\sum_{j=\ell}^m j d_j e^{jw}\right)^2$ is
    $$
    2\sum_{t=1}^{p}t(m-t)d_td_{m-t}=8m^2a^2\left(\frac{H_2}{H_1}\right)^{m}\sum_{t=1}^{p}
    \dbinom{\frac{1}{2}}{m-t}\dbinom{\frac{1}{2}}{t}.
    $$
In addition, \cite[Page~84]{Gould} gives
    $$
    0=\sum_{t=0}^{2p+1}\dbinom{\frac{1}{2}}{m-t}\dbinom{\frac{1}{2}}{t}
    =2\dbinom{\frac{1}{2}}{2p+1} +\sum_{t=1}^{p}\dbinom{\frac{1}{2}}{m-t}\dbinom{\frac{1}{2}}{t}
   +\sum_{t=p+1}^{2p}\dbinom{\frac{1}{2}}{m-t}\dbinom{\frac{1}{2}}{t}.
    $$
It yields that
    $$
    \sum_{t=1}^{p}
    \dbinom{\frac{1}{2}}{m-t}\dbinom{\frac{1}{2}}{t}
    =-\dbinom{\frac{1}{2}}{2p+1} =-\dbinom{\frac{1}{2}}{m}.
    $$

Comparing the coefficients of $e^{sw}$ in \eqref{y''.eq}, we conclude that
    $$
    m^2d_m+2dmd_m=\dbinom{\frac{1}{2}}{m}8m^2a^2\left(\frac{H_2}{H_1}\right)^{m},
    $$
where $d=\deg{\Psi}$. Thus,
    $$
    d_m\left(\frac{H_2}{H_1}\right)^{m}=\frac{m+2d}{2m\dbinom{\frac{1}{2}}{m}}\in\mathbb{Q}.
    $$
Let $x=H_1/H_2e^w=e^w/\tau$, and $g(x)=y(w(x))$. Assume that $g(0)=\Psi(0)=1$. It is obvious that $g(x)$ is an exponential polynomial of order $s=m$, and $g$ is a solution of
    \begin{equation}\label{g.eq}
    x^2g''(x)+x g'(x)-m^2(d_m\tau^m)^2(x^{2m}+x^{2m-1})g(x)=0,
    \end{equation}
where $g'(x)=\frac{dg}{dx}$. Substituting $g(x)=\sum_{n=0}^\infty \gamma_nx^n$ into \eqref{g.eq}, we have $\gamma_0=g(0)=1$, $\gamma_k=0$ for $1\leq k\leq 2m-2$,
    $$
    \gamma_{2m-1}=\frac{m^2d_m^2\tau^{2m}}{(2m-1)^2},
    $$
and
    $$
    \gamma_{k}=\frac{m^2d_m^2\tau^{2m}}{k^2}\left(\gamma_{k-2m}+\gamma_{k-2m+1}\right)\quad\text{for}~k\geq 2m.
    $$
Hence, it yields that $\gamma_k$ is a non-negative rational number for any $k\in\N^+$, and $g(x)\to+\infty$ as $x\to+\infty$.

Assume that $m$ is even. Then $d_m\tau^m<0$. It follows that $g(x)\to 0$ as $x\to+\infty$.
It is contradiction. The only possibility is $m$ is odd. Hence, $d_m\tau^m>0$, and $g(x)\to 0$ as $x\to-\infty$. Choose $x=re^{i\pi/m}$, and put $h(r)=g(re^{i\pi/m})$. It yields by \eqref{g.eq} that
    \begin{equation}\label{r.eq}
    r^2h''(r)+rh'(r)=m^2(d_m\tau^m)^2(r^{2m}+r^{2m-1}e^{-i\pi/m})h(r),
    \end{equation}
where $h'(r)=\frac{dh}{dr}$. Since the Taylor coefficients of $g(x)$ are all real, it gives
    $$
    \overline{h(r)}=\overline{g(re^{i\pi/m})}=g(re^{-i\pi/m}).
    $$
We see that $\overline{h(r)}\to 0$ as $r\to\infty$. Multiplying both sides of equation \eqref{r.eq} by $\overline{h(r)}/r$,
we integrate with respect to $r$ over $(0,\infty)$. Integrating by parts, we have
    $$
    -\int_0^\infty r|h'(r)|^2\,dr=m^2(d_m\tau^m)^2\int_0^\infty(r^{2m-1}+r^{2m-2}e^{-i\pi/m})|h(r)|^2\,dr.
    $$
Taking the imaginary part of both sides of the equation above, it follows that
    $$
    0=-m^2(d_m\tau^m)^2\sin\left(\frac{\pi}{m}\right)\int_0^\infty r^{2m-2}|h(r)|^2\,dr,
    $$
which is impossible.

Therefore, we prove that $\lambda(f)=\infty$.

\appendix

\section{Appendix.}

We recall the definition and basic property of $R$-sets, which are used as exceptional sets in angular domain estimates throughout this paper.

\begin{definition}
A set $E\subset\mathbb{C}$ is called an $R$-set if $E$ can be expressed as a union of a sequence of closed discs $\{D(z_n, r_n)\}_{n=1}^{\infty}$ such that
\[
|z_n|\to\infty,\qquad \sum_{n=1}^{\infty} r_n<\infty.
\]
\end{definition}

It follows immediately that every $R$-set is also an $\mathcal{E}$-set. Recall that an $\mathcal{E}$-set is a union of discs $\{D(z_n,r_n)\}$ with $|z_n|\to\infty$ and $\sum_{n=1}^{\infty} r_n/|z_n|<\infty$. Since
\[
\sum_{n=1}^{\infty}\frac{r_n}{|z_n|}\le \sum_{n=1}^{\infty} r_n
\]
for all sufficiently large $n$, the inclusion $R$-set is $\mathcal{E}$-set holds, whereas the converse is false in general.

The fundamental angular property we use is as follows: the set of arguments $\theta\in[0,2\pi)$ for which the ray $\arg z=\theta$ intersects infinitely many discs of an $R$-set has linear Lebesgue measure zero. That is, almost every ray emanating from the origin meets only finitely many discs of $R$-set.

\begin{lemma}\label{connect.lemma}
Let $G=\{z:\alpha<\arg z<\beta\}$ be an open angular sector with $0<\beta-\alpha<2\pi$, and let $E_1$ be an $R$-set. Then there exists an $R$-set $E_2$ with $E_1\subset E_2$ such that the set $G\setminus E_2$ is connected.
\end{lemma}

\begin{proof}
Write
\[
E_1=\bigcup_{n=1}^{\infty}D(a_n,r_n),
\qquad |a_n|\to\infty,
\qquad \sum_{n=1}^{\infty}r_n<\infty .
\]
Put
\[
\eta=\sin\frac{\beta-\alpha}{2}>0.
\]
We use the following elementary observation: if two closed discs of radii \(r\) and \(s\) intersect, then their union is contained in a disc of radius at most \(r+s\). Hence, by repeatedly merging intersecting or tangent discs, we obtain discs \(D(b_j,s_j)\) such that
\[
E_1\subset \bigcup_jD(b_j,s_j),\qquad
\sum_j s_j<\infty,
\]
and the resulting discs have pairwise disjoint closures.

Since \(s_j\to0\) and \(|b_j|\to\infty\), we have
\[
\frac{s_j}{|b_j|}\longrightarrow0.
\]
After absorbing finitely many exceptional discs, together with all discs meeting a sufficiently large disc centered at the origin, into one disc \(D(0,R_0)\), we may assume that
\[
E_2=D(0,R_0)\cup\bigcup_jD(b_j,s_j)
\]
is an \(R\)-set, \(E_1\subset E_2\), the remaining discs are mutually disjoint and disjoint from \(D(0,R_0)\), and
\[
s_j<\eta |b_j| .
\]
A disc \(D(b_j,s_j)\) which meets both boundary rays
\[
\arg z=\alpha,\qquad \arg z=\beta
\]
must satisfy
\[
s_j\ge |b_j|\sin\frac{\beta-\alpha}{2}
=\eta |b_j|,
\]
which is impossible. Thus no remaining disc joins the two boundary rays of \(G\).

Now \(G\setminus D(0,R_0)\) is connected. Given \(z,w\in G\setminus E_2\), choose a polygonal path in \(G\setminus D(0,R_0)\) joining them. This path meets only finitely many of the remaining discs. Since these discs are mutually disjoint and none joins the two boundary rays, each intersection can be bypassed along the boundary of the corresponding disc. We therefore obtain a path from \(z\) to \(w\) contained in \(G\setminus E_2\).
\end{proof}

\medskip
\noindent
\emph{Z.-T. Wen}\\
\textsc{Shantou University, Department of Mathematics,\\
Daxue Road No.~243, Shantou 515063, China}\\
\texttt{e-mail: zhtwen@stu.edu.cn}

\medskip
\noindent
\emph{S.-J.~Wu}\\
\textsc{Shantou University, Department of Mathematics,\\
Daxue Road No.~243, Shantou 515063, China}\\
\texttt{e-mail: 20sjwu@stu.edu.cn}


\footnotesize
\begin{thebibliography}{99}


\bibitem{BankLangley1987}
S. Bank and J. Langley,
\emph{On the oscillation of solutions of certain linear differential equations in the complex domain},
Proc. Edinb. Math. Soc., II. Ser. \textbf{30} (1987), 455-469.

\bibitem{BankLaine1982}
S. Bank and I. Laine,
\emph{On the oscillation theory of $f''+Af=0$ where $A$ is entire},
Trans. Amer. Math. Soc. \textbf{273} (1982), no.~1, 351--363.

\bibitem{BankLaine}
S. Bank and I. Laine,
\emph{Representations of solutions of periodic second order linear
differential equations},
J. Reine Angew. Math. \textbf{344} (1983), 1--21.

\bibitem{BankLaineLangley1986}
S. Bank, I. Laine and J. Langley,
On the frequency of zeros of solutions of second order linear differential equations,
Result. Math. \textbf{10} (1986), 8--24 .

\bibitem{BankLaineLangley1989}
S. Bank, I. Laine and J. Langley,
\emph{Oscillation results for solutions of linear differential
equations in the complex domain},
Results Math. \textbf{16} (1989), no.~1--2, 3--15.

\bibitem{BergweilerEremenko2017}
W. Bergweiler and A. Eremenko,
\emph{On the Bank--Laine conjecture},
J. Eur. Math. Soc. \textbf{19} (2017), no.~6, 1899--1909.


\bibitem{ChiangLaineWang1995}
Y.-M. Chiang, I. Laine and S. Wang,
\emph{Oscillation results for some linear differential equations},
Math. Scand. \textbf{77} (1995), no.~2, 209--224.

\bibitem{Chen2003}
 Z. Chen,
 \emph{On the hyper order of solutions of higher order differential equations},
Chinese Ann. Math. Ser. B \textbf{24} (2003), no. 4, 501--508.

%

\bibitem{Gould}
H. W. Gould,
\emph{Some generalizations of Vandermonde's convolution},
Amer. Math. Monthly \textbf{63} (1956), 84--91.

\bibitem{Gundersen}
G. G. Gundersen,
\emph{Estimates for the logarithmic derivative of a meromorphic
function, plus similar estimates},
J. London Math. Soc. (2) \textbf{37} (1988), 88--104.

\bibitem{GH2026}
G. Gundersen and J. Heittokangas,
\emph{Entire Functions: Classical and Modern Results, Special Functions and Complex Differential Equations},
De Gruyter, Berlin,  2026.

\bibitem{HITW}
J. Heittokangas, K. Ishizaki, K. Tohge and Z.-T. Wen,
\emph{Value distribution of exponential polynomials and their role in
the theories of complex differential equations and oscillation
theory},
Bull. Lond. Math. Soc. \textbf{55} (2023), no.~1, 1--77.

\bibitem{HILT2019}
J. Heittokangas, K. Ishizaki, I. Laine and K. Tohge,
\emph{Complex oscillation and nonoscillation results},
Trans. Amer. Math. Soc. \textbf{372} (2019), no.~9, 6161--6182.

\bibitem{HILT2021}
J. Heittokangas, K. Ishizaki, I. Laine and K. Tohge,
\emph{Exponential polynomials in the oscillation theory},
J. Differential Equations \textbf{272} (2021), 911--937.

\bibitem{Ishizaki1997}
K. Ishizaki,
\emph{An oscillation result for a certain linear differential equation
of second order},
Hokkaido Math. J. \textbf{26} (1997), no.~2, 421--434.

\bibitem{IshizakiTohge1997}
K. Ishizaki and K. Tohge,
\emph{On the complex oscillation of some linear differential
equations},
J. Math. Anal. Appl. \textbf{206} (1997), no.~2, 503--517.

\bibitem{Won}
K.-H. Kwon,
\emph{On the growth of entire functions satisfying second order linear differential equations}.
Bull. Korean Math. Soc. \textbf{33} (1996), no. 3, 487--496.

\bibitem{Laine}
I. Laine,
\emph{Nevanlinna Theory and Complex Differential Equations},
de Gruyter Studies in Mathematics, Vol.~15,
Walter de Gruyter, Berlin, 1993.
%
%

\bibitem{Langley}
J. K. Langley,
\emph{On complex oscillation and a problem of Ozawa},
Kodai Math. J. \textbf{9} (1986), 430--439.



\bibitem{Wu}
S.-J. Wu,
\emph{Complex oscillation of solutions of second order linear differential equations},
submitted.

\bibitem{YangYi}
C.-C. Yang and H.-X. Yi,
\emph{Uniqueness Theory of Meromorphic Functions},
Mathematics and its Applications, Vol.~557,
Kluwer Academic Publishers, Dordrecht, 2003.

\bibitem{Zhang}
Y. Zhang,
\emph{On the oscillation of certain second-order linear differential
equations},
Proc. Roy. Soc. Edinburgh Sect. A \textbf{153} (2023), 1965--1992.

\bibitem{ZhangCMFT}
Y. Zhang,
\emph{On entire function  $e^{p(z)}\int^z_0\beta(t)e^{-p(t)}dt$  with applications to Tumura-Clunie equations and complex dynamics},
Comput. Methods Funct. Theory
\textbf{23} (2023), no. 2, 213--235.
\end{thebibliography}
\end{document}